\documentclass[11pt]{article}
\usepackage{amssymb,amsthm,amsmath}
\def\qed{\hfill $\Box$}

\newcommand\pf{\smallbreak\noindent \texttt{Proof}. }

\begin{document}

\newtheorem{thm}{Theorem}[section]
\newtheorem{prop}[thm]{Proposition}
\newtheorem{lem}[thm]{Lemma}
\newtheorem{cor}[thm]{Corollary}
\newtheorem{ex}[thm]{Example}
\renewcommand{\thefootnote}{*}

\title{\bf On some relationships between the centers and the derived ideal in Leibniz 3-algebras}

\author{\textbf{P.~Ye.~Minaiev, O.~O.~Pypka}\\
Oles Honchar Dnipro National University, Dnipro, Ukraine\\
{\small e-mail: minaevp9595@gmail.com, sasha.pypka@gmail.com}}
\date{}

\maketitle

\begin{abstract}
One of the classic results of group theory is the so-called Schur theorem. It states that if the central factor-group $G/\zeta(G)$ of a group $G$ is finite, then its derived subgroup $[G,G]$ is also finite. This result has numerous generalizations and modifications in group theory. At the same time, similar investigations were conducted in other algebraic structures. In 2016, L.A. Kurdachenko, J. Otal and O.O. Pypka proved an analogue of Schur theorem for Leibniz algebras: if central factor-algebra $L/\zeta(L)$ of Leibniz algebra $L$ has finite dimension, then its derived ideal $[L,L]$ is also finite-dimensional. Moreover, they also proved a slightly modified analogue of Schur theorem: if the codimensions of the left $\zeta^{l}(L)$ and right $\zeta^{r}(L)$ centers of Leibniz algebra $L$ are finite, then its derived ideal $[L,L]$ is also finite-dimensional. One of the generalizations of Leibniz algebras is the so-called Leibniz $n$-algebras. Therefore, the question of proving analogs of the above results for this type of algebras naturally arises. In this article, we prove the analogues of the two mentioned theorems for Leibniz 3-algebras.
\end{abstract}



\thispagestyle{empty}

\section{Introduction.}
Let $L$ be an algebra over a field $F$ with the binary operations $+$ and $[-,-]$. Then $L$ is called a \textit{left Leibniz algebra} if it satisfies the \textit{left Leibniz identity}:
$$[a,[b,c]]=[[a,b],c]+[b,[a,c]]$$
for all $a,b,c\in L$. Leibniz algebras first appeared in the paper \cite{BA1965}, but the term ``Leibniz algebra'' appears in the book \cite{L1992} and article \cite{L1993}. In \cite{LP1993}, the authors conducted an in-depth study on Leibniz algebras properties. The theory of Leibniz algebras has developed very intensely in many different directions.

Note that Lie algebras present a partial case of Leibniz algebras. Conversely, if $L$ is a Leibniz algebra in which $[a,a]=0$ for every element $a\in L$, then it is a Lie algebra. Thus, Lie algebras can be characterized as anticommutative Leibniz algebras.

One of the key tendencies in the development of Leibniz algebra theory is the search for analogues of the basic results of Lie algebra theory. At the same time, there are very significant differences between these two types of algebras (see, for example, \cite{CPSY2019}).

We recall some necessary definitions. Let $L$ be a Leibniz algebra over a field $F$. If $A,B$ are subspaces of $L$, then $[A,B]$ will denote a subspace generated by all elements $[a,b]$ where $a\in A$, $b\in B$.

A subspace $A$ of $L$ is called a \textit{subalgebra} of $L$, if $[a,b]\in A$ for every $a,b\in A$. It follows that $[A,A]\leq A$. 
A subalgebra $A$ of $L$ is called a \textit{left} (respectively \textit{right}) \textit{ideal} of $L$, if $[b,a]\in A$ (respectively $[a,b]\in A$) for every $a\in A$, $b\in L$. In other words, if $A$ is a left (respectively right) ideal of $L$, then $[L,A]\leq A$ (respectively $[A,L]\leq A$). A subalgebra $A$ of $L$ is called an \textit{ideal} of $L$ (more precisely, \textit{two-sided ideal}) if it is both a left ideal and a right ideal. If $A$ is an ideal of $L$, we can consider a factor-algebra $L/A$. It is not hard to see that this factor-algebra also is a Leibniz algebra.

The \textit{left} $\zeta^{l}(L)$ and \textit{right} $\zeta^{r}(L)$ \textit{centers} of a Leibniz algebra $L$ are defined by the rules:
\begin{align*}
\zeta^{l}(L)&=\{a\in L|\ [a,b]=0\ \mbox{for all }b\in L\},\\
\zeta^{r}(L)&=\{a\in L|\ [b,a]=0\ \mbox{for all }b\in L\}.
\end{align*}
It is not hard to prove that the left center of $L$ is an ideal, but it is not true for the right center. The right center is a subalgebra of $L$, and, in general, the left and right centers are different. They even may have different dimensions (see, for example, \cite{KKPS2017}).

The \textit{center} $\zeta(L)$ of $L$ is the intersection of the left and right centers, that is
$$\zeta(L)=\{a\in L|\ [a,b]=0=[b,a]\ \mbox{for all }b\in L\}.$$
Clearly, the center $\zeta(L)$ is an ideal of $L$. In particular, we can consider the factor-algebra $L/\zeta(L)$.

There is a very close connection between the center $\zeta(L)$ of the Leibniz algebra $L$ and its derived ideal $[L,L]$. In \cite{KOP2016A}, the authors prove that if the central factor-algebra $L/\zeta(L)$ has finite dimension $d$, then the derived ideal $[L,L]$ is also finite-dimensional and $dim_{F}([L,L])\leq d^{2}$. Moreover, they proved some modification of this result: if $codim_{F}(\zeta^{l}(L))=d$ and $codim_{F}(\zeta^{r}(L))=r$ are finite, then $dim_{F}([L,L])\leq d(d+r)$.

The first of the mentioned results is a direct analogue of the so-called \cite{KS2016} Schur theorem. More precisely, Schur theorem states that if the central factor-group $G/\zeta(G)$ of a group $G$ is finite, then the derived subgroup $[G,G]$ is also finite. In this formulation, for the first time it appears in the paper \cite{NB1951}. This theorem was obtained also in \cite{BR1952}. This theorem has numerous generalizations and modifications in group theory (see, for, example, \cite{DKP2014,DKP2015,KOP2015,KOP2016B,KP2014,P2017}). Furthermore, similar investigations were conducted in other algebraic structures of a different nature. Thus, analogues of Schur theorem were obtained for modules (see, for example, \cite{KSC2015}), linear groups~\cite{DKO2013}, topological groups~\cite{U1964}, $n$-groups~\cite{G2006}, associative algebras~\cite{RS2017}, Lie algebras~\cite{KPS2015,S1974}, Lie $n$-algebras~\cite{SV2014}, Lie rings~\cite{KPS2018}, Poisson algebras~\cite{KPS2021}.

In \cite{CLP2002}, the authors introduced the concept of Leibniz $n$-algebras. Let $L$ be an $n$-algebra over a field $F$ with the binary operations $+$ and an $n$-linear bracket $[-,\ldots,-]$. Then $L$ is called a \textit{left Leibniz $n$-algebra} if it satisfies the following \textit{left Leibniz $n$-identity}:
$$[b_{1},\ldots,b_{n-1},[a_{1},\ldots,a_{n}]]=\sum_{i=1}^{n}[a_{1},\ldots,a_{i-1},[b_{1},\ldots,b_{n-1},a_{i}],a_{i+1},\ldots,a_{n}]$$
for any $a_{1},\ldots,a_{n},b_{1},\ldots,b_{n-1}\in L$. The theory of Leibniz $n$-algebras is much less developed than the theory of Leibniz algebras. Our goal is to establish connections between these types of algebras. In particular, in this paper we will prove analogs of the two above-mentioned results from the theory of Leibniz algebras for Leibniz 3-algebras.

\section{Preliminary results.}
Let $L$ be a 3-algebra over a field $F$ with the binary operation $+$ and ternary operation $[-,-,-]$. Then $L$ is called a \textit{Leibniz $3$-algebra} (more precisely, a \textit{left Leibniz $3$-algebra}) if for all elements $a,b,c,x,y\in L$ it satisfies the \textit{left Leibniz $3$-identity}:
$$[x,y,[a,b,c]]=[[x,y,a],b,c]+[a,[x,y,b],c]+[a,b,[x,y,c]].$$

If $A,B,C$ are subspaces of $L$, then $[A,B,C]$ will denote a subspace generated by all elements $[a,b,c]$ where $a\in A$, $b\in B$, $c\in C$. As usual, a subspace $A$ of $L$ is called a \textit{subalgebra} of $L$ if $[a,b,c]\in A$ for all elements 
$a,b,c\in A$. In other words, $[A,A,A]\leq A$. A subalgebra $A$ is called a \textit{left} (respectively \textit{middle}, \textit{right}) \textit{ideal} of $L$ if $[x,y,a]\in A$ (respectively $[x,a,y]\in A$, $[a,x,y]\in A$) for every $x,y\in L$, $a\in A$. In other words, if $A$ is a left (respectively middle, right) ideal of $L$, then $[L,L,A]\leq A$ (respectively $[L,A,L]\leq A$, $[A,L,L]\leq A$). A subalgebra $A$ of $L$ is called an \textit{ideal} of $L$ (or \textit{three-sided ideal}) if it is a left, middle and a right ideal of $L$. If $A$ is an ideal of $L$, we can say about the factor-algebra $L/A$. It is not hard to see that $L/A$ is also a Leibniz 3-algebra.

Let $L$ be a Leibniz 3-algebra over a field $F$, $M$ be a non-empty subset of $L$, and $H$ be a subalgebra of $L$. Put
\begin{align*}
Ann_{H}^{l}(M)&=\{a\in H|\ [a,M,M]=\langle0\rangle\},\\
Ann_{H}^{m}(M)&=\{a\in H|\ [M,a,M]=\langle0\rangle\},\\
Ann_{H}^{r}(M)&=\{a\in H|\ [M,M,a]=\langle0\rangle\}.
\end{align*}
The subset $Ann_{H}^{l}(M)$ (respectively $Ann_{H}^{m}(M)$, $Ann_{H}^{r}(M)$) is called the \textit{left} (respectively \textit{middle}, \textit{right}) \textit{annihilator} of $M$ in $H$. The intersection
$$Ann_{H}(M)=Ann_{H}^{l}(M)\cap Ann_{H}^{m}(M)\cap Ann_{H}^{r}(M)$$
is called \textit{annihilator} of $M$ in $H$. We note the following basic properties of annihilators.

\begin{lem}\label{L1}
Let $L$ be a Leibniz $3$-algebra over a field $F$, $M$ be an ideal of $L$, and $H$ be a subalgebra of $L$. Then $Ann_{H}^{l}(M),Ann_{H}^{m}(M),Ann_{H}^{r}(M)$ are subalgebras of $L$.
\end{lem}
\pf Let $A=Ann_{H}^{l}(M)$, $a,b,c\in A$, $x,y\in M$. We have:
$$[[a,b,c],x,y]=[a,b,[c,x,y]]-[c,[a,b,x],y]-[c,x,[a,b,y]].$$
Since $c\in A$, $[c,x,y]=0$. The fact that $M$ is an ideal of $L$ implies that $[a,b,x]\in M$ and $[a,b,y]\in M$, so that $[c,[a,b,x],y]=[c,x,[a,b,y]]=0$. Thus, $[[a,b,c],x,y]=0$.

Now, let $B=Ann_{H}^{m}(M)$, $a,b,c\in B$, $x,y\in M$. We have:
$$[x,[a,b,c],y]=[a,b,[x,c,y]]-[[a,b,x],c,y]-[x,c,[a,b,y]].$$

Since $c\in B$, $[x,c,y]=0$. Furthermore, $M$ is an ideal of $L$, so that $[a,b,x]\in M$ and $[a,b,y]\in M$. Therefore, $[[a,b,x],c,y]=[x,c,[a,b,y]]=0$. Thus,
$$[x,[a,b,c],y]=0.$$

Finally, let $C=Ann_{H}^{r}(M)$, $a,b,c\in C$, $x,y\in M$. We have:
$$[x,y,[a,b,c]]=[[x,y,a],b,c]+[a,[x,y,b],c]+[a,b,[x,y,c]].$$
Since $a,b,c\in C$, $[x,y,a]=[x,y,b]=[x,y,c]=0$. Thus, $[x,y,[a,b,c]]=0$.

Since $Ann_{H}(M)=A\cap B\cap C$, $Ann_{H}(M)$ is also a subalgebra of $L$. \qed

Let $L$ be a Leibniz 3-algebra over a field $F$. Put
$$\zeta(L)=\{a\in L|\ [a,x,y]=[x,a,y]=[x,y,a]=0\ \mbox{for all }x,y\in L\}.$$
The subset $\zeta(L)$ is called the \textit{center} of $L$. We note that the center of $L$ is the annihilator of $L$ in $L$. The left, middle and right annihilators lead us to the following subsets. Put
\begin{align*}
\zeta^{l}(L)&=\{a\in L|\ [a,x,y]=0\ \mbox{for all }x,y\in L\},\\
\zeta^{m}(L)&=\{a\in L|\ [x,a,y]=0\ \mbox{for all }x,y\in L\},\\
\zeta^{r}(L)&=\{a\in L|\ [x,y,a]=0\ \mbox{for all }x,y\in L\}.
\end{align*}
The subset $\zeta^{l}(L)$ (respectively $\zeta^{m}(L)$, $\zeta^{r}(L)$) is called the \textit{left} (respectively \textit{middle}, \textit{right}) \textit{center} of $L$. Obviously,
$$\zeta(L)=\zeta^{l}(L)\cap\zeta^{m}(L)\cap\zeta^{r}(L).$$
In particular, we can define the \textit{lm-center} of $L$ by the rule
$$\zeta^{lm}(L)=\zeta^{l}(L)\cap\zeta^{m}(L).$$

We note that, in general, $\zeta^{l}(L),\zeta^{m}(L),\zeta^{r}(L)$ are not ideals of $L$. However, we have the following result.

\begin{lem}\label{L2}
Let $L$ be a Leibniz $3$-algebra over a field $F$. Then the following assertions hold.
\begin{enumerate}
\item[\upshape(i)] $\zeta^{l}(L),\zeta^{m}(L),\zeta^{r}(L),\zeta(L)$ are subalgebras of $L$.
\item[\upshape(ii)] $\zeta^{lm}(L),\zeta(L)$ are ideals of $L$.
\end{enumerate}
\end{lem}
\pf (i) Since
$$Ann_{L}^{l}(L)=\zeta^{l}(L),Ann_{L}^{m}(L)=\zeta^{m}(L),Ann_{L}^{r}(L)=\zeta^{r}(L),Ann_{L}(L)=\zeta(L),$$
Lemma~\ref{L1} shows that $\zeta^{l}(L),\zeta^{m}(L),\zeta^{r}(L),\zeta(L)$ are subalgebras of $L$.

(ii) Let $D=\zeta^{lm}(L)$, $v\in D$, $a,b,x,y\in L$. We have:
\begin{align*}
[[x,y,v],a,b]&=[x,y,[v,a,b]]-[v,[x,y,a],b]-[v,a,[x,y,b]],\\
[a,[x,y,v],b]&=[x,y,[a,v,b]]-[[x,y,a],v,b]-[a,v,[x,y,b]].
\end{align*}
Since $v\in D$,
\begin{align*}
[[x,y,v],a,b]&=[x,y,0]-0-0=0,\\
[a,[x,y,v],b]&=[x,y,0]-0-0=0.
\end{align*}
This means that $D$ is a left ideal of $L$. Furthermore,
\begin{align*}
[[x,v,y],a,b]&=[x,v,[y,a,b]]-[y,[x,v,a],b]-[y,a,[x,v,b]]\\
&=0-[y,0,b]-[y,a,0]=0,\\
[a,[x,v,y],b]&=[x,v,[a,y,b]]-[[x,v,a],y,b]-[a,y,[x,v,b]]\\
&=0-[0,y,b]-[a,y,0]=0,
\end{align*}
which shows that $D$ is a middle ideal of $L$. Finally, we have:
\begin{align*}
[[v,x,y],a,b]&=[v,x,[y,a,b]]-[y,[v,x,a],b]-[y,a,[v,x,b]]\\
&=0-[y,0,b]-[y,a,0]=0,\\
[a,[v,x,y],b]&=[v,x,[a,y,b]]-[[v,x,a],y,b]-[a,y,[v,x,b]]\\
&=0-[0,y,b]-[a,y,0]=0.
\end{align*}
In other words, $D$ is a right ideal of $L$. Thus, $D$ is an ideal of $L$.

In the same way we can show that the center $\zeta(L)$ of $L$ is an ideal of $L$. \qed

Let $L$ be a Leibniz 3-algebra over a field $F$. A linear transformation $f$ of $L$ is called a \textit{derivation} of $L$ if
$$f([a,b,c])=[f(a),b,c]+[a,f(b),c]+[a,b,f(c)]$$
for all $a,b,c\in L$.

Consider the mapping $l_{a,b}:L\rightarrow L$ defined by the rule $l_{a,b}(x)=[a,b,x]$, $a,b\in L$. We note some basic properties of $l_{a,b}$. Let $x,y,z\in L$, $\lambda\in F$. Then
\begin{align*}
l_{a,b}(x+y)&=[a,b,x+y]=[a,b,x]+[a,b,y]=l_{a,b}(x)+l_{a,b}(y),\\
l_{a,b}(\lambda x)&=[a,b,\lambda x]=\lambda[a,b,x]=\lambda l_{a,b}(x),\\
l_{a,b}([x,y,z])&=[a,b,[x,y,z]]\\
&=[[a,b,x],y,z]+[x,[a,b,y],z]+[x,y,[a,b,z]]\\
&=[l_{a,b}(x),y,z]+[x,l_{a,b}(y),z]+[x,y,l_{a,b}(z)].
\end{align*}
These equalities show that $l_{a,b}$ is a derivation of $L$.

\section{Main result.}
We now present the main result of this paper.
\begin{thm}\label{T1}
Let $L$ be a Leibniz $3$-algebra over a field $F$. If $codim_{F}(\zeta^{lm}(L))=d$ and $codim_{F}(\zeta^{r}(L))=r$ are finite, then $dim_{F}([L,L,L])\leq d^{2}(d+r)$.
\end{thm}
\pf We have: $L=\zeta^{lm}(L)\oplus E$ for some subspace $E$. Choose in $E$ a basis $\{e_{1},\ldots,e_{d}\}$. If $x,y,z$ are arbitrary elements of $L$, then
\begin{align*}
x&=\alpha_{1}e_{1}+\ldots+\alpha_{d}e_{d}+s_{1},\\
y&=\beta_{1}e_{1}+\ldots+\beta_{d}e_{d}+s_{2},\\
z&=\gamma_{1}e_{1}+\ldots+\gamma_{d}e_{d}+s_{3}
\end{align*}
for $\alpha_{i},\beta_{i},\gamma_{i}\in F$, $i\in\{1,\ldots,d\}$, and $s_{j}\in\zeta^{lm}(L)$, $j\in\{1,2,3\}$. Then we have
\begin{align*}
[x,y,z]&=\left[\sum_{1\leq i\leq d}\alpha_{i}e_{i}+s_{1},\sum_{1\leq j\leq d}\beta_{j}e_{j}+s_{2},\sum_{1\leq k\leq d}\gamma_{k}e_{k}+s_{3}\right]\\
&=\sum_{1\leq i,j,k\leq d}\alpha_{i}\beta_{j}\gamma_{k}[e_{i},e_{j},e_{k}]+\sum_{1\leq i,j\leq d}\alpha_{i}\beta_{j}[e_{i},e_{j},s_{3}].
\end{align*}
It follows that the subspace $S$ generated by the elements $[e_{i},e_{j},e_{k}]$, $1\leq i,j,k\leq d$, and the subspaces $[e_{i},e_{j},\zeta^{lm}(L)]$, $1\leq i,j\leq d$, include $[L,L,L]$.

Put $Z=\zeta^{lm}(L)$ and let $a,b$ be arbitrary elements of $L$. We define a mapping $l_{a,b}(x):Z\rightarrow Z$ by the rule $l_{a,b}(z)=[a,b,z]$, $z\in Z$. As we noted above, this mapping is linear, $Im(l_{a,b})=[a,b,Z]$ and $Ker(l_{a,b})=Ann_{Z}^{r}((a,b))$ where $(a,b)$ is an ordered pair of elements $a$ and $b$. Hence
$$[a,b,Z]=Im(l_{a,b})\cong_{F}Z/Ker(l_{a,b})=Z/Ann_{Z}^{r}((a,b)).$$
Since $\zeta^{r}(L)\leq Ann_{L}^{r}((a,b))$, we have $codim_{F}(Ann_{L}^{r}((a,b)))\leq r$ and so
$$dim_{F}([a,b,Z])=dim_{F}(Z/Ann_{Z}^{r}((a,b)))\leq r.$$
In particular, $dim_{F}([e_{i},e_{j},Z])\leq r$ for every $1\leq i,j\leq d$. It follows that the subspace $S$ has dimension at most $d^{3}+d^{2}r=d^{2}(d+r)$, as required. \qed

\begin{cor}\label{C1}
Let $L$ be a Leibniz $3$-algebra over a field $F$. If $codim_{F}(\zeta(L))=d$ is finite, then $dim_{F}([L,L,L])\leq d^{3}$.
\end{cor}

Indeed, we have $\zeta(L)=\zeta^{lm}(L)\cap\zeta^{r}(L)$. Since $codim_{F}(\zeta^{lm}(L))$ and $codim_{F}(\zeta^{r}(L))$ are finite, $codim_{F}(\zeta(L))$ is finite. Then it suffices to take the proof of Theorem~\ref{T1} into account to obtain Corollary~\ref{C1}.

Recall \cite{FV1985} that a 3-algebra $L$ over a field $F$ ($char(F)\neq2$) with the binary operations $+$ and a 3-linear bracket $[-,-,-]$ is called a \textit{Lie $3$-algebra} if it satisfies the following conditions:

(i) Lie 3-bracket is antisymmetric, that is $[a_{1},a_{2},a_{3}]=sign(\sigma)[a_{\sigma(1)},a_{\sigma(2)},a_{\sigma(3)}]$;

(ii) Lie 3-bracket satisfies the generalized Jacobi identity (or Jacobi 3-identity), that is
$$[x,y,[a,b,c]]=[[x,y,a],b,c]+[a,[x,y,b],c]+[a,b,[x,y,c]]$$
for any $a,b,c,x,y\in L$ and any permutation $\sigma\in S_{3}$. Note that since $char(F)\neq2$ the first condition equivalent to $[a_{1},a_{2},a_{3}]=0$ whenever $a_{i}=a_{j}$ for some $i\neq j$, $1\leq i,j\leq3$. Obviously, every Lie 3-algebra is a Leibniz 
3-algebra. Therefore, we can apply the previous results to this partial case. We have the following

\begin{cor}\label{C2}
Let $L$ be a Lie $3$-algebra over a field $F$. If $codim_{F}(\zeta(L))=d$ is finite, then $dim_{F}([L,L,L])\leq\frac{d(d-1)(d-2)}{6}$.
\end{cor}

\end{document}